\documentclass[a4paper,12pt]{amsart}
\usepackage{amssymb}
\usepackage{ifthen}
\usepackage{graphicx}
\usepackage{color,xcolor}
\nonstopmode
\numberwithin{equation}{section}
\newtheorem{thm}{Theorem}
\newtheorem{cor}{Corollary}
\newtheorem{lem}{Lemma}
\newtheorem{prop}{Proposition}

\newtheorem{conj}{Conjecture}
\newtheorem{prob}{Problem}

\theoremstyle{definition}
\newtheorem{defn}{Definition}

\newtheorem{rem}{Remark}

\newenvironment{pf}[1][]{%
 \vskip 1mm
 \noindent
 \ifthenelse{\equal{#1}{}}%
  {{\slshape Proof. }}%
  {{\slshape #1.} }%
 }%
{\qed\medskip}

\newcounter{alphabet}
\newcounter{tmp}
\newenvironment{Thm}[1][]{\refstepcounter{alphabet}%
\bigskip%
\noindent%
{\bf Theorem \Alph{alphabet}}%
\ifthenelse{\equal{#1}{}}{}{ (#1)}%
{\bf .} \itshape}{\vskip 8pt}

\makeatletter
\newcommand{\Ref}[1]{\@ifundefined{r@#1}{}{\setcounter{tmp}{\ref{#1}}\Alph{tmp}}}
\makeatother
\newenvironment{Lem}[1][]{\refstepcounter{alphabet}%
\bigskip%
\noindent%
{\bf Lemma \Alph{alphabet}}%
{\bf .} \itshape}{\vskip 8pt}

\newenvironment{Conjs}[1][]{\refstepcounter{alphabet}%
\bigskip%
\noindent%
{\bf Conjectures \Alph{alphabet}}%
{\bf .} \itshape}{\vskip 8pt}

\newcounter{alphabet2}

\newcommand{\ID}{{\mathbb D}}

\newcommand{\dist}{{\operatorname{dist}}}

\def\be{\begin{equation}}
\def\ee{\end{equation}}

\newcommand{\ben}{\begin{enumerate}}
\newcommand{\een}{\end{enumerate}}

\newcommand{\blem}{\begin{lem}}
\newcommand{\elem}{\end{lem}}
\newcommand{\bthm}{\begin{thm}}
\newcommand{\ethm}{\end{thm}}
\newcommand{\bcor}{\begin{cor}}
\newcommand{\ecor}{\end{cor}}
\newcommand{\beg}{\begin{exam}}
\newcommand{\eeg}{\end{exam}}
\newcommand{\begs}{\begin{examples}}
\newcommand{\eegs}{\end{examples}}
\newcommand{\bdefe}{\begin{defn}}
\newcommand{\edefe}{\end{defn}}
\newcommand{\bprob}{\begin{prob}}
\newcommand{\eprob}{\end{prob}}
\newcommand{\bques}{\begin{ques}}
\newcommand{\eques}{\end{ques}}
\newcommand{\bei}{\begin{itemize}}
\newcommand{\eei}{\end{itemize}}
\newcommand{\bcon}{\begin{conj}}
\newcommand{\econ}{\end{conj}}
\newcommand{\bop}{\begin{op}}
\newcommand{\eop}{\end{op}}

\newcommand{\bas}{\begin{assertion}}
\newcommand{\eas}{\end{assertion}}

\newcommand{\bfa}{\begin{fact}}
\newcommand{\efa}{\end{fact}}

\newcommand{\bca}{\begin{ca}}
\newcommand{\eca}{\end{ca}}

\newcommand{\bst}{\begin{step}}
\newcommand{\est}{\end{step}}

\newcommand{\bsca}{\begin{sca}}
\newcommand{\esca}{\end{sca}}

\newcommand{\bcl}{\begin{cl}}
\newcommand{\ecl}{\end{cl}}

\newcommand{\bmlem}{\begin{mlem}}
\newcommand{\emlem}{\end{mlem}}

\newcommand{\bscl}{\begin{scl}}
\newcommand{\escl}{\end{scl}}

\newcommand{\bcons}{\begin{conjs}}
\newcommand{\econs}{\end{conjs}}

\newcommand{\bprop}{\begin{prop}}
\newcommand{\eprop}{\end{prop}}

\newcommand{\br}{\begin{rem}}
\newcommand{\er}{\end{rem}}
\newcommand{\brs}{\begin{rems}}
\newcommand{\ers}{\end{rems}}
\newcommand{\bo}{\begin{obser}}
\newcommand{\eo}{\end{obser}}
\newcommand{\bos}{\begin{obsers}}
\newcommand{\eos}{\end{obsers}}
\newcommand{\bpf}{\begin{pf}}
\newcommand{\epf}{\end{pf}}
\newcommand{\ba}{\begin{array}}
\newcommand{\ea}{\end{array}}
\newcommand{\beq}{\begin{eqnarray}}
\newcommand{\beqq}{\begin{eqnarray*}}
\newcommand{\eeq}{\end{eqnarray}}
\newcommand{\eeqq}{\end{eqnarray*}}

\newcommand{\ra}{\to}

\newcommand{\ds}{\displaystyle}

\newcounter{minutes}
\divide\time by 60
\newcounter{hours}
\multiply\time by 60 \addtocounter{minutes}{-\time}
\begin{document}

\bibliographystyle{amsplain}
\title []
{Bohr's phenomenon for the classes of Quasi-subordination and $K$-quasiregular harmonic mappings}

%%%%%%%% BEGIN TIMESTAMP
\def\thefootnote{}
\footnotetext{ \texttt{\tiny File:~\jobname .tex,
          printed: \number\day-\number\month-\number\year,
          \thehours.\ifnum\theminutes<10{0}\fi\theminutes}
} \makeatletter\def\thefootnote{\@arabic\c@footnote}\makeatother
%%%%%%%% END TIMESTAMP

\author{Ming-Sheng Liu}
 \address{M.S. Liu, School of Mathematical Sciences, South China Normal University, Guangzhou 510631, China.} \email{liumsh65@163.com}

\author{Saminathan Ponnusamy %${}^{~\mathbf{*}}$
}
\address{S. Ponnusamy, Department of Mathematics,
Indian Institute of Technology Madras, Chennai-600 036, India. }
\email{samy@iitm.ac.in}

\author{Jun Wang}
\address{J. Wang, School of Mathematical Sciences, Fudan University, Shanghai 200433, China.}
\email{majwang@fudan.edu.cn}

%DEPARTMENT OF MATHEMATICS, SHANTOU UNIVERSITY, SHANTOU, GUANGDONG
%515063, PEOPLE¡¯S REPUBLIC OF CHINA E-mail address:
%xtwang@stu.edu.cn (X. Wang)

\keywords{Bohr radius, analytic functions, Harmonic mappings, Convex functions, Subordination,
Quasi-subordination,  $K$-quasiregular mappings\\
%${}^{\mathbf{*}}$ Corresponding author
}

%\thanks{This research is partly supported by Guangdong Natural Science Foundations (Grant No. 2018A030313508).}
\begin{abstract}
In this paper, we investigate the Bohr radius for $K$-quasiregular sense-preserving harmonic mappings $f=h+\overline{g}$
in the unit disk $\mathbb{D}$ such that the translated analytic part $h(z)-h(0)$ is quasi-subordinate to some analytic function. The main aim of this article is to extend and to establish sharp versions of four recent theorems by Liu and Ponnusamy \cite{LP2019}
%[Z.H. Liu and S. Ponnusamy, Bohr Radius for Subordination and $K$-quasiregular Harmonic Mappings,
%Bull. Malays. Math. Sci. Soc., \textbf{42}(2019), 2151--2168.]
and, in particular, we settle affirmatively the two conjectures proposed by them. Furthermore, we establish two refined versions of Bohr's
inequalities and determine the Bohr radius for the derivatives of analytic functions associated with quasi-subordination.
\end{abstract}

\maketitle \pagestyle{myheadings}
\markboth{M-S Liu, S. Ponnusamy, and J. Wang}{Quasi-subordination and $K$-quasiregular harmonic mappings}

\section{Introduction}
Throughout we let $\mathcal B$ denote the class of all analytic functions $\omega$ in the unit disk $\mathbb{D}=\{z\in \mathbb{C}:|z|<1\}$
such that $|\omega(z)|\leq 1$ for all $z\in \mathbb{D}$, and let ${\mathcal B}_0=\{\omega\in {\mathcal B}:\, \omega(0)=0\}$.
%Functions in ${\mathcal B}_0$ are often referred as Schwarz functions.
Bohr's inequality says that if $f\in {\mathcal B}$
and $f(z)=\sum_{n=0}^{\infty}a_nz^n$, then, for the majorant series $M_f(r)=\sum_{k=0}^{\infty}|a_k|r^k$ of $f$, we have
\begin{equation*}
M_{f_0}(r)=\sum_{k=1}^{\infty}|a_k|r^k\leq 1-|a_0| =\dist (f(0), \partial \ID)
\end{equation*}
for all $z\in \mathbb{D}$ with $|z|=r\leq\frac{1}{3}$, where $f_0(z)=f(z)-f(0)$. This inequality was discovered by Bohr in 1914 \cite{AMK1914}., Bohr actually obtained the inequality for $|z|\leq\frac{1}{6}$. Later, F. Wiener, M. Riesz and I. Schur, independently established the inequality for $|z|\leq\frac{1}{3}$ and showed that $1/3$ is sharp. See \cite{AAP2016}, \cite[Chapter 8]{GarMasRoss-2018}  and the references therein.
%\cite{PPS2002,S1927,T1962}.
Few other proofs are also available in the literature. %\cite{PS2004,PS2006}.
The Bohr radius has been discussed  for certain power series in $\mathbb{D}$, as well as for analytic functions from $\mathbb{D}$ into specific
domains, such as convex domains, simply connected domains, the punctured unit disk, the exterior of the closed unit disk,
and concave wedge-domains. The analogous Bohr radius has also been studied for harmonic and starlike log-harmonic mappings in $\mathbb{D}$.
In particular, in \cite{KayPon1}, the authors settled the conjecture of  Ali et al. \cite{AliBarSoly} about the Bohr radius for
odd functions from $\mathcal B$.  In the year 2000, powered Bohr inequality was initiated by Djakov and Ramanujan  \cite{DjaRaman-2000}
and a conjecture related to their work was settled in \cite{KayPon_AAA18} affirmatively.
Bohr's idea naturally extends to functions of several complex variables \cite{A2000,A2005,AAP2016,BK1997,DjaRaman-2000}. Several other
aspects of Bohr's inequality may be obtained from \cite{Aizen-12,AAP2016,AKP2019,BDK-2004,BombBour-196,BombBour-2004,DFOOS,DIN,GarMasRoss-2018,KM,LSX2018,PPS2002,PS2004,PS2006} and the references therein.

In this article, we shall consider Bohr's radius for complex-valued $K$-quasiregular harmonic mappings of the unit disk $\mathbb{D}$. In order to state our main results, we recall the following notions and notations (see \cite{AKP2019, R1970}).

\bdefe
Let $f(z)$ and $g(z)$ be analytic in $\mathbb{D}$. We say that
\begin{enumerate}
\item $f(z)$ is subordinate to $g(z)$ in $\mathbb{D}$, written by $f(z)\prec g(z)$ or $f\prec g$,
if there exists an $\omega\in {\mathcal B}_0$  such that $f(z)=g(\omega (z))$ for $z\in\ID$.
Furthermore, if $g(z)$ is univalent in $\mathbb{D}$, then we have the following relation
$$ f(z)\prec g(z)\Longleftrightarrow f(0)=g(0),~f(\mathbb{D})\subset g(\mathbb{D}).
$$
\item  $f(z)$ is majorized by $g(z)$ in $\mathbb{D}$, denoted by $f(z)\ll g(z)$ or $f\ll g$,
if $|f(z)|\leq |g(z)|$ for $z \in\mathbb{D}$.
\item $f$ is quasi-subordinate to $g$ (relative to $\Phi$), denoted by $f(z)\prec_q g(z)$ in $\mathbb{D}$,
if there exists a $\Phi\in {\mathcal B}$ and an $\omega\in {\mathcal B}_0$ such that
%\begin{equation}
$f(z) = \Phi(z)g(\omega(z))
%\label{liu11} \end{equation}
$ for $z\in\ID$.
\end{enumerate}
\edefe

Evidently if either $f\prec g$ or $|f(z)|\leq |g(z)|$ in $\mathbb{D}$, then $f(z)\prec_q g(z)$ in $\mathbb{D}$. Thus, the notion of quasi-subordination generalizes both the concept of subordination and the principle of majorization.

A harmonic mapping in $\mathbb{D}$ is a complex-valued function $f=u+iv$ of $z=x+iy$ in $\mathbb{D}$, which satisfies the
Laplace equation $\triangle f=4f_{z\overline{z}}=0$.
%, where $f_{z}=\frac{1}{2}(f_{x}-if_{y})$, $f_{\overline{z}}=\frac{1}{2}(f_{x}+if_{y})$.
It follows that every harmonic mapping $f$ admits a representation of the form $f=h+\overline{g}$, where $h$ and $g$ are analytic in $D$.
This representation is unique up to an additive constant. It is convenient to assume that $f(0)=g(0)$.
The Jacobian $J_{f}$ of $f$ is given by $J_{f}(z)=|h'(z)|^{2}-|g'(z)|^{2}$.

We say that $f$ is sense-preserving in $\mathbb{D}$ if $J_{f}(z)>0$ in $\mathbb{D}$. Consequently,
a harmonic mapping $f$ is locally univalent and sense-preserving in $\mathbb{D}$ if and only if $J_{f}(z)>0$ in $\mathbb{D}$; or equivalently if $h'\neq 0$ in $\mathbb{D}$ and the dilatation $\omega_{f}:=\frac{g'}{h'}$ of $f$ has the property that $|\omega_{f}|<1$ in $\mathbb{D}$ \cite{L1936}.

%In order to state the results about Bohr radius for quasi-conformal harmonic mappings, we need to introduce some notations. For harmonic mappings $f$ in $D$, we use the following standard notations:
%\begin{equation*}  \Lambda_{f}(z)=\max_{0\leq\theta\leq2\pi}|f_{z}(z)+e^{-2i\theta}f_{\overline{z}}(z)|=|f_{z}(z)|+|f_{\overline{z}}(z)|\end{equation*}
%and\begin{equation*}  \lambda_{f}(z)=\min_{0\leq\theta\leq2\pi}|f_{z}(z)+e^{-2i\theta}f_{\overline{z}}(z)|=\bigg||f_{z}(z)|-|f_{\overline{z}}(z)|\bigg|,\end{equation*}
%so that if $f$ is locally univalent and sense-preserving, then
%\begin{equation*}  J_{f}=\lambda_{f}\Lambda_{f}=|f_{z}|^{2}-|f_{\overline{z}}|^{2}>0.\end{equation*}
%A sense-preserving homeomorphism $f$ from the unit disk $D$ onto $\Omega$, contained in the Sobolev class $W_{loc}^{1,2}(D)$, is said to be a $K$-quasiregular mapping if, for $z\in D$,
%\begin{equation*}  \frac{|f_{z}|+|f_{\overline{z}}|}{|f_{z}|-|f_{\overline{z}}|}=\frac{1+|\omega_{f}(z)|}{1-|\omega_{f}(z)|}\leq K,~~ie.,~~|\omega_{f}(z)|\leq k=\frac{K-1}{K+1},\end{equation*}
%where $K\geq1$ so that $k\in[0,1)$\cite{LV1973,V1988}. %The majorant series is defined by $M_{f}(r)=\sum_{n=0}^{\infty}|a_n|r^n$.

If a locally univalent and sense-preserving harmonic mapping $f=h+\overline{g}$ satisfies the condition
$ |\omega_{f}(z)| \leq k<1 ,
$
then $f$ is called $K$-quasiregular harmonic mapping on $\mathbb{D}$, where $K=\frac{1+k}{1-k}\geq 1$
(cf. \cite{Kalaj2008,Martio1968}). %, and also \cite{Zhu2016} for some recent investigation on harmonic $K$-quasiregular self-mapping of $\ID$).
Obviously $k\rightarrow 1$ corresponds to the limiting case $K\rightarrow \infty$. Note that when $k=1$, the condition on
the dilatation of $f$ becomes
$|\omega _f(z)|\leq 1$ in which case the Jacobian could be zero at some point. Thus, it is worth pointing out that
our results below cover this case as well as the case where $f$ is sense-preserving. A harmonic extension of the classical
Bohr theorem was established in \cite{AM2010,KP2017,Kayumov2}.

\bdefe
We say that $f=h+\overline{g}\in H_{K,\,h\prec_q \varphi}(\ID)$ if it is a $K$-quasiregular sense-preserving harmonic mapping of $\mathbb{D}$
and has the power series form
$$f(z)=h(z)+\overline{g(z)}=\sum_{n=0}^{\infty} a_{n}z^{n}+\overline{\sum_{n=1}^{\infty }b_{n}z^{n}}, \quad z\in\ID,
$$
together with an additional condition that $h\prec_q \varphi$,  where $k=\frac{K-1}{K+1}$. If $h\prec \varphi$, then
we simply write $H_{K,\,h\prec \varphi}(\ID)$ instead of $H_{K,\,h\prec_q \varphi}(\ID)$ by suppressing the subscript $`q$'.

Similarly, we can define $H_{K,\,h\ll \varphi}(\ID)$ by replacing quasi-subordination condition $h\prec_q \varphi$ by
the majorization condition $h\ll \varphi$.
\edefe

%{\color{red} \br
%In above definition, {\color{red} $K$-quasiregular} harmonic mapping need not be a homeomorphism,
%which differs from the usual one in the literature.
%\er}

Recently, Liu and Ponnusamy \cite{LP2019} have investigated the class $H_{K,\,h\prec \varphi}(\ID)$ and
obtained the following results.

\begin{Thm} {\rm (\cite{LP2019})}\label{ThA}
%Suppose that $f(z)=h(z)+\overline{g(z)}=\sum_{n=0}^{\infty} a_{n}z^{n}+\overline{\sum_{n=1}^{\infty }b_{n}z^{n}}$ is a
$K$-quasiregular sense-preserving harmonic mapping in $\mathbb{D}$, $k=\frac{K-1}{K+1}$ and $h\prec \varphi$.
For $f=h+\overline{g}\in H_{K,\,h\prec \varphi}(\ID)$ with $k=(K-1)/(K+1)$, we have the following:
\begin{enumerate}
\item If $\varphi$ is analytic and univalent in $\mathbb{D}$, then
\begin{eqnarray}\label{liu11a}
\sum_{n=1}^{\infty}(|a_{n}|+|b_{n}|)r^{n}\leq {\rm dist}\, (\varphi(0), \partial\varphi (\mathbb{D}))
\end{eqnarray}
for $|z|=r\leq r_u$, where $r_u=r_u(k)$ is the root of the equation
$$
(1-r)^2-4r(1+k\sqrt{1+r})=0
$$
in the interval $(0, 1)$.

\item If $\varphi$ is univalent and convex in $\mathbb{D}$, then \eqref{liu11a} holds for
$|z| = r \leq\frac{K + 1}{5K + 1}$. The result is sharp.

\item If  $\varphi$ is univalent and convex in $\mathbb{D}$ and $b_1=g'(0)=0$, then \eqref{liu11a} holds for
 $|z|=r\leq r_c$, where $r_c=r_c(k)$ is the positive root of the equation
\begin{equation*}\label{eq:r}
\frac{r}{1-r}+\frac{kr^2}{1-r^2}\sqrt{\left(\frac{1+r^2}{1-r^2}\right)\left(\frac{\pi^2}{6}-1\right)}
=\frac{1}{2}.
\end{equation*}
The number $r_c(k)$ cannot be replaced by a number larger than $\rho:=\rho_c(k)$, where $\rho$ is the positive root of the equation
\begin{equation}\label{eq:R}
\frac{2(1+k)\rho}{1-\rho}+2k\ln (1-\rho)=1.
\end{equation}

\item If $\varphi$ is analytic and univalent in $\mathbb{D}$, $h(0)=0$, and $b_1=g'(0)=0$, then  \eqref{liu11a} holds
for $|z|=r\leq r_s$, where $r_s=r_s(k)$ is the positive real root of the equation
\begin{equation*}\label{rootS}
\frac{r}{(1-r)^2}+\frac{kr^2}{(1-r^2)^2}\sqrt{\left(\frac{r^6+11 r^4+11 r^2+1}{1-r^2}\right)\left(\frac{\pi^2}{6}-1\right)}=\frac{1}{4}
\end{equation*}
%\begin{equation}\label{rootS}
%(1-r)^2- 4 r \left(1+\sqrt{r (1+r)}\right)=0
%\end{equation}
in the interval $(0,1)$. The number $r_s(k)$ cannot be replaced by a number larger than $\rho =\rho_s(k)$,
where $\rho$ is the positive root of the equation
\begin{equation}\label{eq:R1}
\frac{\rho(1-k+2k\rho)}{(1-\rho)^2}-k\ln  (1-\rho)=\frac{1}{4}.
\end{equation}
\end{enumerate}
\end{Thm}

One of the important special cases is when $K\ra \infty$, i.e. $k\ra 1$. Thus, the authors in \cite{LP2019}
proposed the following two conjectures.

%{\bf Conjecture E}(\cite{LP2019})\quad
\begin{Conjs}\label{conj1}%\begin{conj}\label{conj1}
{\it
Suppose that $f(z)=h(z)+\overline{g(z)}=\sum_{n=0}^{\infty}a_n z^n+\overline{\sum_{n=2}^{\infty}b_n z^n}$ is a
%locally univalent and
sense-preserving harmonic mapping in $\ID$ and $h\prec \varphi$.
\begin{enumerate}
\item[{\rm (a)}]
If $\varphi$ is univalent and convex in $\ID$, then
\begin{equation}\label{conj-eq1}
\sum_{n=1}^{\infty}|a_n|r^n+\sum_{n=2}^{\infty}|b_n|r^n\leq \dist (\varphi(0),\partial\varphi(\ID))
\end{equation}
for $|z|=r\leq \rho_c= 0.299823\cdots$, where $\rho_c$ is the positive root of the equation
\begin{eqnarray*}
\frac{4r}{1-r}+2\ln (1-r)=1,%\label{liu14}
\end{eqnarray*}(compare with \eqref{eq:R} with $k=1$).

\item[{\rm (b)}] If $\varphi$ is univalent in $\ID$, then the inequality \eqref{conj-eq1} holds
%\begin{equation*}
%\sum_{n=1}^{\infty}|a_n|r^n+\sum_{n=2}^{\infty}|b_n|r^n\leq \dist (0,\partial\varphi(\ID))
%\end{equation*}
for $|z|=r\leq \rho_s = 0.161353\cdots$, where $\rho_s$ is the positive real root of the equation
\begin{eqnarray*}
\frac{2r^2}{(1-r)^2}-\ln (1-r)=\frac{1}{4}
%\label{liu15}
\end{eqnarray*}
(compare with \eqref{eq:R1} with $k=1$).
\end{enumerate}
}
%\end{conj}Under the hypotheses of parts {\rm (3)} and {\rm (4)} of Theorem \Ref{ThA} with $k=1$, $\rho_c=0.299823\cdots$ and  $\rho_s=0.161353\cdots$ are sharp, where $\rho_c$ and $\rho_s$ are respectively the  positive root of the equations
%\begin{eqnarray*} \frac{4r}{1-r}+2\ln (1-r)=1,\label{liu14} \end{eqnarray*}(compare with \eqref{eq:R} with $k=1$) and
%\begin{eqnarray*}\frac{2r^2}{(1-r)^2}-\ln (1-r)=\frac{1}{4}\label{liu15}\end{eqnarray*}(compare with \eqref{eq:R1} with $k=1$).
\end{Conjs}
One of the aims of this article is to prove sharp versions of Theorem \Ref{ThA}(3) and (4) which in turn imply that
Conjectures \Ref{conj1} are true. In fact, we prove these in a general setting along with the sharp version of Theorem \Ref{ThA}(1).

Inspired by Theorem \Ref{ThA} and the notion of quasi-subordination in the setting of Bohr's inequality, as discussed in \cite{AKP2019},
we obtain the following results. % after thinking and theoretical research. Next we will prove the following theorems:

\bthm\label{LSW1-th1}
Let $f=h+\overline{g}\in H_{K,\,h_0\prec_q \varphi _0}(\ID)$, $k=\frac{K-1}{K+1}$, $h_0(z)=h(z)-h(0)$ and
$\varphi_0(z)=\varphi(z)-\varphi (0)$. We have the following:
\begin{enumerate}
\item If $\varphi$ is analytic and univalent in $\mathbb{D}$, then
\be\label{LSW1-eq2}
\sum_{n=1}^{\infty}(|a_{n}|+|b_{n}|)r^{n}\leq {\rm dist}\, (\varphi(0), \partial\varphi (\mathbb{D}))
\ee
for $|z|=r\leq r_u$, where $r_u=r_u(k)=\frac{1}{2k+3+\sqrt{(2k+3)^2-1}}\in (0,1/3)$. The result is sharp.

\item If $\varphi$ is univalent and convex in $\mathbb{D}$, then \eqref{LSW1-eq2} holds for
$|z| = r \leq\frac{K + 1}{5K + 1}$.  The result is sharp.
\end{enumerate}
\ethm

\br
%\begin{enumerate}
%\item
The conclusion of Theorem \ref{LSW1-th1}(1) continues to hold under the assumption that $f=h+\overline{g}\in H_{K,\,h\prec \varphi }(\ID)$
and thus, Theorem \ref{LSW1-th1} contains a sharp version of Theorem \Ref{ThA}(1) or \cite[Theorem 3]{LP2019}, namely,
with the subordination $h_0\prec \varphi _0$ (which is equivalent to $h\prec \varphi $) in place of $h_0\prec_q \varphi _0$.
Note that in the case of quasi-subordination,  $h_0\prec_q \varphi _0$ is not equivalent to $h\prec_q \varphi$ unless $h(0)=\varphi (0)=0$.
In particular, if we set $k=0$ (i.e. $K=1$) in the case of subordination,  then we get the result of Abu-Muhanna \cite[Theorem 1]{AM2010} as
a special case of Theorem \ref{LSW1-th1}.
%\end{enumerate}
\er

\br
Recall that the notion of quasi-subordination generalizes both the concept of subordination and the principle of majorization
and thus, Theorem \ref{LSW1-th1}(2) is an extension of Theorem \Ref{ThA}(2)  or  \cite[Theorem 1]{LP2019}. In particular,
the conclusion of Theorem \ref{LSW1-th1}(2) holds if we replace the condition $h_0\prec_q \varphi _0$ by the
majorization condition $|h_0(z)|\leq |\varphi_0(z)|$ for $z\in\mathbb{D}$.
\er

\bthm\label{LSW1-th2}
Let $f=h+\overline{g}\in H_{K,\,h\prec \varphi}(\ID)$, $b_1=g'(0)=0$, and $k=\frac{K-1}{K+1}$. We have the following:
\begin{enumerate}
\item If $\varphi$ is univalent and convex in $\mathbb{D}$, then
\be\label{LSW1-eq3}
\sum_{n=1}^{\infty}|a_{n}| r^{n}+\sum_{n=2}^{\infty}|b_{n}|r^{n}\leq {\rm dist}\, (\varphi(0), \partial\varphi (\mathbb{D}))
%\quad\mbox{ for }\, |z| = r \leq\frac{K + 1}{5K + 1}.
\ee
for $|z|=r\leq r_c$, where $r_c=r_c(k)$ is the unique positive root in $(0, 1/3)$ of the equation
\be\label{liu16}
\frac{2(1+k)r}{1-r}+2k\ln (1-r)=1.
\ee
The result is sharp.

\item If $\varphi$ is analytic and univalent in $\mathbb{D}$, %and $h(0)=0=\varphi (0)$,
then \eqref{LSW1-eq3} holds for $|z|=r\leq r_s$, where $r_s=r_s(k)$ is the unique positive root in $(0, 1/3)$ of the equation
\be\label{liu17}
\frac{r(1-k+2k r)}{(1-r)^2}-k\ln (1-r)=\frac{1}{4}.
\ee
The result is sharp.
\end{enumerate}
\ethm

\br
Theorem \ref{LSW1-th2}(1) is the sharp version of Theorem \Ref{ThA}(3)  or  \cite[Theorem 2]{LP2019}.
Also, Theorem \ref{LSW1-th2}(2) is the sharp version of Theorem \Ref{ThA}(4) or \cite[Theorem 4]{LP2019}.
Setting $k=0$ in Theorem \ref{LSW1-th2}(1), we also get the classical version of the Bohr theorem.
Finally, we remark that the whole proof of Theorem \ref{LSW1-th2} can be imitated to establish Conjectures \Ref{conj1} only by replacing $k$ with $1$. %the case $k=1$ of Theorem \ref{LSW1-th2} settles both Conjectures \Ref{conj1} affirmatively.
\er

The paper is organized as follows. In Section \ref{LPW1-sec3}, we present the proof of Theorems \ref{LSW1-th1} and \ref{LSW1-th2}.
In Section \ref{LSW1-sec4}, we state and prove two theorems which extend two recent results of
Ponnusamy et al. \cite{PVW2019,PVW201911} from the case of analytic functions to the case of sense-preserving harmonic mappings.
Finally, in Section \ref{LPW1-sec5}, we investigate the Bohr radius for the derivatives of analytic functions in the setting
of quasi-subordination.

\section{Preliminaries}%\label{LPW1-sec1}
In order to establish %the proofs of
Theorems \ref{LSW1-th1} and \ref{LSW1-th2}, we need the following lemmas.
It is easy to obtain the following two well-known lemmas from the latest monograph of
 Avkhadiev and  Wirths \cite{AvWir-09}. See also \cite[pp. 195--196]{D1983} and \cite{AM2010,K2003}.

\begin{Lem} \label{LPW1-LemC}
Let $\varphi$ be an analytic univalent map from $\mathbb{D}$ onto a simply connected domain
$\Omega =\varphi(\mathbb{D})$. We have the following:
\begin{enumerate}
\item $\ds
\frac{1}{4}|\varphi'(0)|\leq {\rm dist}\, (\varphi(0), \partial\Omega)\leq |\varphi'(0)|.$
\item If $g(z) = \sum_{n=0}^{\infty} b_{n}z^{n}\prec\varphi(z)$, then
$\ds |b_n|\leq n|\varphi'(0)|\leq 4n\, {\rm dist}\, (\varphi(0), \partial\Omega).
$
\end{enumerate}
\end{Lem}

\begin{Lem} \label{LPW1-LemD}
Let $\varphi$ be an analytic univalent map from $\mathbb{D}$ onto a convex domain
$\Omega =\varphi(\mathbb{D})$. We have the following:
\begin{enumerate}
\item $\ds
\frac{1}{2}|\varphi'(0)|\leq {\rm dist}\, (\varphi(0), \partial\Omega)\leq |\varphi'(0)|.
$
\item
If $g(z) = \sum_{n=0}^{\infty} b_{n}z^{n}\prec\varphi(z)$, then
$\ds |b_n|\leq |\varphi'(0)|\leq 2{\rm dist}\, (\varphi(0), \partial\Omega).
$
\end{enumerate}
\end{Lem}

The following two lemmas will play a key role in the proofs of our main results in
Section \ref{LPW1-sec3}.

\begin{Lem} \label{LPW1-LemE}
 {\rm (Alkhaleefah et al. \cite{AKP2019})}
Let $f(z)$ and $g(z)$ be two analytic functions in $\mathbb{D}$ with the Taylor
series expansions $f(z) = \sum_{n=0}^{\infty} a_{n}z^{n}$ and $g(z) = \sum_{n=0}^{\infty} b_{n}z^{n}$ for $z\in\mathbb{D}$.
We have the following:
\begin{enumerate}
\item[(1)] If $f(z)\prec_q g(z)$ in $\mathbb{D}$, then
\begin{eqnarray*}
\sum_{n=0}^{\infty}|a_{n}|r^{n}\leq \sum_{n=0}^{\infty}|b_{n}|r^{n}\quad\mbox{ for all }\, r \leq\frac{1}{3}.
\end{eqnarray*}
\item[(2)] If  $|g'(z)|\leq k|h'(z)|$ in $\mathbb{D}$ for some $k \in (0, 1]$, then
\begin{eqnarray*}
\sum_{n=1}^{\infty}|b_{n}|r^{n}\leq k\sum_{n=1}^{\infty}|a_{n}|r^{n}\quad\mbox{ for all }\, r \leq\frac{1}{3}.
\end{eqnarray*}
\end{enumerate}
\end{Lem}
\bpf
The proof of the first part of Lemma \Ref{LPW1-LemE} is available %well-known
 in  \cite{AKP2019}  while the second part follows easily from this.
Indeed, by assumption, we obtain that $g'(z)\prec_q k h'(z)$ which quickly gives from Lemma \Ref{LPW1-LemE}(1) that
\begin{eqnarray*}
\sum_{n=1}^{\infty}n|b_{n}|r^{n-1}\leq \sum_{n=1}^{\infty}k n\, |a_{n}|r^{n-1}\quad\mbox{ for all }\, r \leq\frac{1}{3}
\end{eqnarray*}
and integrating this with respect to $r$ gives the desired inequality.
\epf

%\begin{Lem} \label{LPW1-lemF}
%{\color{red} There exists a unique positive root $r_c=r_c(k), r_s=r_s(k)$ in $(0, 1/3)$ of the equation (\ref{liu16}) and (\ref{liu17}) respectively.}
%\end{Lem}
%\bpf
%{\color{red} Let $F(r)=\frac{2(1+k)r}{1-r}+2k\ln(1-r)-1$, then it is evident that $F(r)$ is continuous in $[0, 1/3]$, and
%\begin{eqnarray*}
%&&F(0)=-1<0,\quad F(\frac{1}{3})=k\ln\frac{4e}{9}>0,\\
%&&F'(r)=\frac{2+2kr}{(1-r)^2}>0,\quad r\in [0,\frac{1}{3}],
%\end{eqnarray*}
%so it follows from the intermediate value theorem that the equation (\ref{liu16}) has a unique root $r_c=r_c(k)$ in $(0, 1/3)$.
%
%Similar, we may verify that the equation (\ref{liu17}) has a unique root $r_s=r_s(k)$ in $(0, 1/3)$.}
%\epf
%\hfill $\Box$

%{\bf Lemma 2.3}(\cite{RP2017}) \, Let $|a|<1$ and $0<R\leq1$. If $g(z)=\sum_{k=0}^{\infty}b_{k}z^{k}$ is analytic and satisfies the inequality $|g(z)|\leq1$ in $D$. Then the following sharp inequality holds:
%\begin{equation}  \sum_{k=1}^{\infty}|b_{k}|^{2}R^{pk}\leq R^{p}\frac{(1-|b_{0}|^{2})^{2}}{1-|b_{0}|^{2}R^{p}}.  \label{liu02}\end{equation}
%{\bf Lemma 2.4}(\cite{GK2003}) \, Suppose that $p(z)=\sum_{n=0}^{\infty}p_nz^n$ is analytic in $D$ satisfy $Re~p(z)>0$ in $D$. Then$|p_{n}|\leq 2{\rm Re}~p_{0}$ for all $n=1,2,\cdots$.

%{\bf Lemma 2.5}(\cite{GK2003}) \, If $h(z)=\sum_{k=0}^{\infty}c_{k}z^{k}$ is analytic in $D$ such that $|h(z)|\leq1$ for all $z\in D$. Then $|c_{k}|\leq1-|c_{0}|^{2}$ for all $n=1,2,\cdots$.

\section{The proofs of Theorems \ref{LSW1-th1} and \ref{LSW1-th2}}\label{LPW1-sec3}

%In this section, we will give the proofs of Theorems $1\sim 4$ mentioned in Section 1.

\subsection{Proof of Theorem  \ref{LSW1-th1}}
Assume that $f=h+\overline{g}\in H_{K,\,h_0\prec_q \varphi _0}(\ID)$. Then
$f(z)=h(z)+\overline{g(z)}=\sum_{n=0}^{\infty} a_{n}z^{n}+\overline{\sum_{n=1}^{\infty }b_{n}z^{n}}$,
where $|g'(z)|\leq k|h'(z)|$ in $\mathbb{D}$ for some $k\in [0,1)$ and $h(z)-h(0)\prec_q \varphi(z)-\varphi (0)$ in $\mathbb{D}$.

(1) By assumption, $\varphi(z)= \sum_{n=0}^{\infty} c_{n}z^{n}$ is analytic and univalent in $\mathbb{D}$. Now,
by Lemmas \Ref{LPW1-LemC} and \Ref{LPW1-LemE}(1), we obtain respectively the inequalities
\begin{eqnarray*}
|c_n|\leq n |\varphi'(0)|\leq 4n\, {\rm dist}\, (\varphi(0), \partial\varphi (\mathbb{D})),\quad\mbox{ for }\,  n=1,2,\ldots ,
\end{eqnarray*}
and
\begin{eqnarray*}
\sum_{n=1}^{\infty}|a_{n}|r^{n}\leq \sum_{n=1}^{\infty}|c_{n}|r^{n}\quad\mbox{ for }\, r \leq\frac{1}{3},
\end{eqnarray*}
which by the previous inequality leads to
\begin{eqnarray}
\sum_{n=1}^{\infty}|a_{n}|r^{n}\leq %4 \dist (\varphi(0), \partial(\varphi (\mathbb{D}))\sum_{n=1}^{\infty}n r^{n}=
4 \dist (\varphi(0), \partial\varphi (\mathbb{D}))\frac{r}{(1-r)^2}\quad\mbox{ for all }\, r \leq\frac{1}{3}.
\label{liu31}
\end{eqnarray}
Again, as $|g'(z)|\leq k|h'(z)|$, Lemma \Ref{LPW1-LemE}(2) and \eqref{liu31} show that
\begin{eqnarray*}
\sum_{n=1}^{\infty}|b_{n}|r^{n}\leq k\sum_{n=1}^{\infty}|a_{n}|r^{n}\leq 4k \dist (\varphi(0), \partial\varphi (\mathbb{D}))\frac{r}{(1-r)^2}\quad\mbox{ for all }\, r \leq\frac{1}{3}.
\end{eqnarray*}
Consequently, by combining the last two inequalities, we have
\begin{eqnarray*}
\sum_{n=1}^{\infty}(|a_{n}|+|b_{n}|)r^{n}\leq 4(1+k) \dist (\varphi(0), \partial\varphi (\mathbb{D}))\frac{r}{(1-r)^2}
\end{eqnarray*}
which is less than or equal to $\dist (\varphi(0), \partial\varphi (\mathbb{D}))$ if and only if
\begin{eqnarray*}
4(1+k)r-(1-r)^2= (r-r_u )(2k+3+\sqrt{(2k+3)^2-1}-r)\leq 0.
\end{eqnarray*}
The last inequality holds for $r\leq r_u$, where $r_u=r_u(k)=\frac{1}{2k+3+\sqrt{(2k+3)^2-1}}\in (0,1/3)$ and $k=\frac{K-1}{K+1}$.% Substituting $k =\frac{K-1}{K+1} $ gives the desired result.

To prove the sharpness, we consider $f=h+\overline{g}$ such that $g'(z)=k\lambda h'(z)$,  where $\lambda\in\mathbb{D}$,
\begin{eqnarray*}
h(z)=a_0+\frac{z}{(1-z)^2} ~\mbox{ and }~ \varphi(z)=c_0+\frac{z}{(1-z)^2}=c_0+\sum_{n=1}^{\infty}n z^{n}.
\end{eqnarray*}
Then $\dist (\varphi(0), \partial\varphi (\mathbb{D}))=1/4$
%$$
%\dist (\varphi(0), \partial\varphi (\mathbb{D}))=\frac{1}{4}, ~\mbox{ $a_n=n$ and $b_n=k\lambda n$ for $n\geq 1$},
%$$
and so it is easy to see that
\begin{eqnarray*}
\sum_{n=1}^{\infty}(|a_{n}|+|b_{n}|)r^{n} =%\sum_{n=1}^{\infty}(1+k|\lambda|)nr^n=
(1+k|\lambda|)\frac{r}{(1-r)^2}
\end{eqnarray*}
which is bigger than or equal to $1/4$ if and only if $4(1+k|\lambda|)r-(1-r)^2\geq 0.$
%\begin{eqnarray*}
%4(1+k|\lambda|)r-(1-r)^2&=& \left (r-\frac{1}{2k|\lambda|+3+\sqrt{(2k|\lambda|+3)^2-1}}\right )\\
%&&\cdot \left (2k|\lambda|+3+\sqrt{(2k|\lambda|+3)^2-1}-r\right )\geq 0.
%\end{eqnarray*}
Solving the last inequality shows that the number $r_u=\frac{1}{2k+3+\sqrt{(2k+3)^2-1}}$ cannot be improved since $|\lambda|$ could be chosen so close to 1 from left. This completes the proof of the first part (1).

(2) For the proof of the second part, we just need to assume that $\varphi$ is convex and then proceed
 with the above method of proof, but using Lemma \Ref{LPW1-LemD} in place of Lemma \Ref{LPW1-LemC}. This change after
 minor computation  leads to
\begin{eqnarray*}
\sum_{n=1}^{\infty}(|a_{n}|+|b_{n}|)r^{n}\leq 2(1+k) \dist (\varphi(0), \partial\varphi (\mathbb{D}))\frac{r}{1-r}%\quad\mbox{ for }\, |z| = r \leq\frac{K + 1}{5K + 1}.
\end{eqnarray*}
which is less than or equal to $\dist (\varphi(0), \partial\varphi (\mathbb{D}))$ for $r\leq\frac{1}{3+2k}\leq\frac{1}{3}$. Substituting $k =\frac{K-1}{K+1} $ gives the desired result.

Again, to prove the sharpness, we consider $f=h+\overline{g}$ such that
\begin{eqnarray*}
h(z)=a_0+\frac{z}{1-z},\quad \varphi(z)=c_0+\frac{z}{1-z}=c_0+\sum_{n=1}^{\infty}z^{n},
\end{eqnarray*}
and $g'(z)=k\lambda h'(z)$, where $\lambda\in\mathbb{D}$. Then a simple computation yields
$$\dist (\varphi(0), \partial\varphi (\mathbb{D}))=\frac{1}{2} ~\mbox{ and }~g(z)=k\lambda\frac{z}{1-z}
$$
so that for this function we have
\begin{eqnarray*}
\sum_{n=1}^{\infty}(|a_{n}|+|b_{n}|)r^{n} %=\sum_{n=1}^{\infty}(1+k|\lambda|)r^n
= (1+k|\lambda|)\frac{r}{1-r}
\end{eqnarray*}
which is bigger than or equal to $1/2$ if and only if
\begin{eqnarray*}
r\geq\frac{1}{3+2k|\lambda|}=\frac{K+1}{3K+3+2|\lambda|(K-1)}.
\end{eqnarray*}
This shows that the number $\frac{K + 1}{5K + 1}$ cannot be improved since $|\lambda|$ could be chosen so
close to $1$ from left. This completes the proof of the second part. \hfill $\Box$

\subsection{Proof of Theorem  \ref{LSW1-th2}}
Suppose that $f=h+\overline{g}\in H_{K,\,h\prec \varphi}(\ID)$. We consider the first part of the theorem,
where $b_1=g'(0)=0$,  $h\prec \varphi$ and $\varphi$ is univalent and convex in $\mathbb{D}$.
It follows from Lemma \Ref{LPW1-LemD}(2) that $|a_n|\leq |\varphi'(0)|$ for $n\ge 1$, and thus
%\begin{eqnarray*}
%|a_n|\leq |\varphi'(0)|\leq 2\dist (\varphi(0), \partial\varphi (\mathbb{D}))\quad\mbox{ for } n\geq 1,
%\end{eqnarray*}
%so that
\begin{eqnarray}
\sum_{n=1}^{\infty}|a_{n}|r^{n}\leq %|\varphi'(0)|\sum_{n=1}^{\infty}r^{n}=
|\varphi'(0)|\frac{r}{1-r}.
\label{liu35}
\end{eqnarray}
Because $g'(0)=0$, by Schwarz's lemma, we obtain that $\omega=\frac{g'}{h'}$ is analytic in $\mathbb{D}$ and $|\omega(z)|\leq k|z|$ in $\mathbb{D}$. Therefore, we have $|g'(z)|\leq |kz h'(z)|$, or $g'(z)\prec_q kz h'(z)$ in $\mathbb{D}$. %where $0\leq k< 1$,

By Lemma \Ref{LPW1-LemE}(1), we have
\begin{eqnarray*}
\sum_{n=1}^{\infty}n|b_{n}|r^{n-1}\leq \sum_{n=1}^{\infty}k n |a_{n}|r^{n}\leq k|\varphi'(0)|\sum_{n=1}^{\infty}nr^{n}=
k|\varphi'(0)|\frac{r}{(1-r)^2}
\quad\mbox{ for  }\, r \leq\frac{1}{3}.
\end{eqnarray*}
Integrate this inequality on $[0, r]$, where $r\leq 1/3$, we obtain
\begin{eqnarray*}
\sum_{n=1}^{\infty}|b_{n}|r^{n}
%&=&\int_0^r\sum_{n=1}^{\infty}n|b_{n}|t^{n-1}dt\\
%&\leq &\int_0^r \sum_{n=1}^{\infty}k n|a_{n}|t^{n}dt=\sum_{n=1}^{\infty}k n|a_{n}|\int_0^r t^{n}dt\\
%&=&k\sum_{n=1}^{\infty}\frac{n}{n+1}|a_{n}|r^{n+1}\\
%&\leq& k|\varphi'(0)|\sum_{n=1}^{\infty}\frac{n}{n+1}r^{n+1}
&\leq& k|\varphi'(0)| \int_0^r \frac{t}{(1-t)^2}\,dt
= k|\varphi'(0)| \left (\ln  (1-r)+\frac{r}{1-r}\right ).
\end{eqnarray*}%for all $r \leq\frac{1}{3}$.
Consequently, by combining (\ref{liu35}) with the last inequality, we find that
\begin{eqnarray*}
\sum_{n=1}^{\infty}|a_{n}| r^{n}+\sum_{ n=1}^{\infty} |b_{n}| r^{n}&\leq & |\varphi'(0)|\left [\frac{r}{1-r}+k\left (\ln (1-r)+\frac{r}{1-r}\right )\right ]\\
&\leq & 2\dist (\varphi(0), \partial\varphi (\mathbb{D})) \left [\frac{(1+k)r}{1-r}+k\ln (1-r)\right ]\\
&\leq & \dist (\varphi(0), \partial\varphi (\mathbb{D})),
\end{eqnarray*}
where the last inequality holds if and only if
\begin{eqnarray*}
\frac{2(1+k)r}{1-r}+2k\ln (1-r)\leq 1.
\end{eqnarray*}
The above inequality holds for $r\leq r_c(k)$, where $r_c(k)$ is  the unique positive root in $(0, 1/3)$ of equation (\ref{liu16}).
% {\color{red} from Lemma \Ref{LPW1-lemF}}. % \Ref{LPW1-LemC}Substituting $k =\frac{K-1}{K+1} $ gives the desired result.
In order to verify the fact, we let
$$F(r)=\frac{2(1+k)r}{1-r}+2k\ln(1-r)-1.
$$
Then $F(r)$ is continuous in $[0, 1/3]$,
$$F(0)=-1<0,~F(1/{3})=k\ln\frac{4e}{9}>0,~ F'(r)=\frac{2+2kr}{(1-r)^2}>0 ~\mbox{ for } r\in [0,{1}/{3}],
$$
and therefore, it follows from the intermediate value theorem that the equation (\ref{liu16}) has a unique
root $r_c=r_c(k)$ in $(0, 1/3)$.

To prove the sharpness, we consider $f=h+\overline{g}$, where  $h, g$ and $\varphi$ are such that  $g'(z)=k z h'(z)$ and
$h(z)=\varphi(z)=1/(1-z)$.
%\begin{eqnarray*}
%h(z)=\varphi(z)=\frac{1}{1-z}=\sum_{n=0}^{\infty} z^{n}.
%\end{eqnarray*}
Then for these choices we find that $ \dist (\varphi(0), \partial\varphi (\mathbb{D}))=1/2$ and it is easy to compute the
corresponding sum
%$$ \dist (\varphi(0), \partial\varphi (\mathbb{D}))=\frac{1}{2},
%$$
%$a_n=1$ for $n\geq 1$, and $b_n=\frac{k(n-1)}{n}$ for $n\geq 2$, so that
\begin{eqnarray*}
\sum_{n=1}^{\infty}|a_{n}| r^{n}+\sum_{n=1}^{\infty} |b_{n}| r^{n}
%=\sum_{n=1}^{\infty}r^{n}+\sum_{n=2}^{\infty}\frac{k(n-1)}{n} r^{n}
=\frac{(1+k)r}{1-r}+k\ln (1-r),
\end{eqnarray*}
which is less than or equal to $1/2$ only for $r\leq r_c(k)$, where $r_c(k)$ is the unique positive root in $(0, 1/3)$ of the equation (\ref{liu16}). This shows that the number $r_c(k)$ cannot be improved. This completes the proof of the first part of the theorem.

Now we consider the second part, where $\varphi$ is analytic and univalent in $\mathbb{D}$. %, $\varphi(0)=h(0)=0$.
It follows from Lemma \Ref{LPW1-LemC} that $|a_n|\leq |\varphi'(0)| n\leq
4 n\, \dist (\varphi (0), \partial\varphi (\mathbb{D}))$
for $n\geq 1$, and thus
\begin{eqnarray}
\sum_{n=1}^{\infty}|a_{n}|r^{n}\leq %|\varphi'(0)|\sum_{n=1}^{\infty}n r^{n}=
|\varphi'(0)|\frac{r}{(1-r)^2}.
\label{liu36}
\end{eqnarray}

As in the proof of the previous case, we have
\begin{eqnarray*}
\sum_{n=1}^{\infty}n|b_{n}|r^{n-1}\leq \sum_{n=1}^{\infty}k n |a_{n}|r^{n}\leq k|\varphi'(0)|\sum_{n=1}^{\infty}n^2r^{n}=
k|\varphi'(0)|\frac{r(1+r)}{(1-r)^3}
\quad\mbox{ for  }\, r \leq\frac{1}{3}
\end{eqnarray*}
and thus, by integration, we obtain easily that
\begin{eqnarray*}
\sum_{n=1}^{\infty}|b_{n}|r^{n}&\leq &
%k\sum_{n=1}^{\infty}\frac{n}{n+1}|a_{n}|r^{n+1}\\
%&\leq& k|\varphi'(0)|\sum_{n=1}^{\infty}\frac{n^2}{n+1}r^{n+1}\\ &=&
k|\varphi'(0)| \left ( \frac{2r^2-r}{(1-r)^2}-\ln (1-r)\right ) \quad\mbox{ for  }\, r \leq\frac{1}{3}.
\end{eqnarray*}
Consequently, by combining (\ref{liu36}) with the last inequality, we find that
\begin{eqnarray*}
\sum_{n=1}^{\infty}|a_{n}| r^{n}+\sum_{n=1}^{\infty} |b_{n}| r^{n}&\leq & |\varphi'(0)|\left [\frac{r}{(1-r)^2}+k\left (\frac{2r^2-r}{(1-r)^2}-\ln (1-r)\right )\right ]\\
&\leq & 4\dist (\varphi(0), \partial\varphi (\mathbb{D})) \left [\frac{r(1-k+2k r)}{(1-r)^2}-k\ln (1-r)\right ]\\
&\leq & \dist (\varphi(0), \partial\varphi (\mathbb{D})),
\end{eqnarray*}
where the last inequality holds if and only if
\begin{eqnarray*}
\frac{r(1-k+2k r)}{(1-r)^2}-k\ln (1-r)\leq \frac{1}{4},
\end{eqnarray*}
which holds for $r\leq r_s(k)$, where $r_s(k)$ is the unique positive root in $(0, 1/3)$ of equation
(\ref{liu17})--a fact which is easy to verify as in the proof of the case (\ref{liu16}) above. % Substituting $k =\frac{K-1}{K+1} $ gives the desired result.

To prove the sharpness, we consider the function $f=h+\overline{g}$, where
\begin{eqnarray*}
h(z)=\varphi(z)=\frac{z}{(1-z)^2}=a_0+\sum_{n=1}^{\infty} n z^{n},
\end{eqnarray*}
and $g'(z)=k z h'(z)$. Then we find that $\dist (\varphi(0), \partial\varphi (\mathbb{D}))=\frac{1}{4}$, and as before we have
%$$|a_n|=n, \mbox{ and } |b_n|=k(n+\frac{1}{n}-2)\, \mbox{ for } n\geq 2, $$
%and thus, we have
\begin{eqnarray*}
\sum_{n=1}^{\infty}|a_{n}| r^{n}+\sum_{n=2}^{\infty} |b_{n}| r^{n}=\sum_{n=1}^{\infty}n r^{n}+ k\sum_{n=2}^{\infty}\Big(n+\frac{1}{n}-2\Big) r^{n}=\frac{r(1-k+2k r)}{(1-r)^2}-k\ln (1-r),
\end{eqnarray*}
which is less than or equal to $1/4$ only in the case where $r\leq r_s(k)$, where $r_s(k)$ is the unique positive root in $(0, 1/3)$ of equation
(\ref{liu17}). This shows that the number $r_s(k)$ cannot be improved. The proof of the theorem  is complete. \hfill $\Box$

\section{Improved Bohr's phenomenon associated with quasi-subordination}\label{LSW1-sec4}

 Recently, Ponnusamy et al. \cite{PVW2019,PVW201911} established several refined versions of Bohr's inequality in the case of bounded analytic functions. In this section, following the ideas of \cite{PVW2019,PVW201911}, we will discuss improved Bohr's phenomenon for two classes of sense-preserving harmonic mappings associated with quasi-subordination.

\bthm\label{LSW1-th3}
Suppose that $f(z)=h(z)+\overline{g(z)}=\sum_{n=0}^{\infty} a_{n}z^{n}+\overline{\sum_{n=1}^{\infty }b_{n}z^{n}}$ is a sense-preserving harmonic mapping in $\mathbb{D}$ and $h(z)-h(0)\prec_q \varphi(z)-\varphi (0)$ in $\mathbb{D}$, where $\varphi(z)$ is univalent and convex in $\mathbb{D}$. Also, let $\lambda=\dist (\varphi(0), \partial\varphi (\mathbb{D}))<1$ and  $\|f_0\|_r=\sum_{n=1}^{\infty}(|a_{n}|^2+|b_{n}|^2)r^{2n}$, where $f_0(z)=f(z)-f(0)$. Then
\begin{eqnarray*}
T_f(r):=\sum_{n=1}^{\infty}(|a_{n}|+|b_{n}|)r^{n}+\left(\frac{1}{2-\lambda}+\frac{r}{1-r}\right )\|f_0\|_r\leq \lambda
~\mbox{ for $|z|=r\leq r_*$,}
\end{eqnarray*}
where $r_*\approx 0.15867508$ is the unique root in $(0, 1)$ of equation
$$5r^3-9r^2-5r+1=0.
$$
Moreover, for any $\lambda\in (0, 1)$, there exists a uniquely defined $r_0\in (r_*, \frac{1}{5})$ such that $T_f(r)\leq \lambda$ for $r\in [0, r_0]$. The radius $r_0$ can be calculated as the solution of the equation
\begin{eqnarray*}
\Phi (\lambda, r)=8r^3\lambda^2-(13r^3+7r^2-5r+1)\lambda +10 r^3-2r^2-10r+2=0.
\end{eqnarray*}
\ethm
\bpf Let $\varphi(z)= \sum_{n=0}^{\infty} c_{n}z^{n}$. Then, as before, Lemma \Ref{LPW1-LemD} implies that
$|c_n|\leq |\varphi'(0)|\leq 2\lambda$ for $n\geq 1$. Because $h(z)-h(0)\prec_q \varphi(z)-\varphi (0)$ and $|g'(z)|\leq |f'(z)|$ for $z\in \mathbb{D}$, we have
\begin{eqnarray*}
\sum_{n=1}^{\infty}|b_{n}|r^{n}\leq \sum_{n=1}^{\infty}|a_{n}|r^{n}\leq \sum_{n=1}^{\infty}|c_{n}|r^{n}\quad \mbox{ for } r\leq \frac{1}{3}.
\end{eqnarray*}
In addition we also have
\begin{eqnarray*}
\sum_{n=1}^{\infty}|a_{n}|^2 r^{2n}\leq \sum_{n=1}^{\infty}|c_{n}|^2r^{2n} =: \|\varphi_0\|_r ~\mbox{ and }~
\sum_{n=1}^{\infty}n^2|b_{n}|^2 r^{2n-2}\leq \sum_{n=1}^{\infty}n^2|a_{n}|^2r^{2n-2}
\end{eqnarray*}
so that
\begin{eqnarray*}
\sum_{n=1}^{\infty}|b_{n}|^2 r^{2n}\leq\sum_{n=1}^{\infty}|a_{n}|^2 r^{2n}\leq \sum_{n=1}^{\infty}|c_{n}|^2r^{2n}=\|\varphi_0\|_r.
\end{eqnarray*}
Consequently,
$\|f_0\|_r\leq 2\|\varphi_0\|_r $ and thus, $T_f(r)\leq 2T_{\varphi}(r)$ for $r\leq 1/3$ only, where
$$T_{\varphi}(r)=\sum_{n=1}^{\infty}|c_{n}|r^{n}+\left (\frac{1}{2-\lambda}+\frac{r}{1-r}\right)\|\varphi_0\|_r.
$$
Clearly, the desired conclusion follows if we can show that $T_{\varphi}(r)\leq \lambda /2$.
Finally, because $|c_n|\leq 2\lambda$ for $n\geq 1$, we have
\begin{eqnarray*}
2T_{\varphi}(r)&\leq &4\lambda\sum_{n=1}^{\infty}r^{n}+\Big(\frac{1}{2-\lambda}+\frac{r}{1-r}\Big)8\lambda^2\sum_{n=1}^{\infty}r^{2n}\\
&= &\frac{4\lambda r}{1-r}+\frac{1+(1-\lambda)r}{(2-\lambda)(1-r)}\cdot \frac{8\lambda^2r^2}{1-r^2}\\
&= &\lambda-\lambda\Big[\frac{1-5r}{1-r}-\frac{8\lambda r^2(1+(1-\lambda)r)}{(2-\lambda)(1-r)(1-r^2)}\Big]\\
&= &\lambda-\lambda\cdot\frac{\Phi (\lambda, r)}{(2-\lambda)(1-r)(1-r^2)},
\end{eqnarray*}
where
\begin{eqnarray*}
\Phi (\lambda, r)&=&(1-5r)(2-\lambda)(1-r^2)- 8\lambda r^2(1+(1-\lambda)r)\\
&=& 8r^3\lambda^2-\left[(1-5r)(1-r^2)+8r^2(1+r)\right]\lambda +2(1-5r)(1-r^2).
\end{eqnarray*}
It is easy to see that $\Phi (\lambda, r)<0$ for $r>1/5$ and $0<\lambda\leq 1$, so that
$$\lambda-\lambda\cdot\frac{\Phi (\lambda, r)}{(2-\lambda)(1-r)(1-r^2)}>\lambda$$
for $r>1/5$ and $0<\lambda\leq 1$.

We claim that $\Phi (\lambda, r)\geq 0$ for every $r\leq r_*$ and for $\lambda\in (0, 1]$. In fact, we have
\begin{eqnarray*}
\frac{\partial^2\Phi (\lambda, r)}{\partial\lambda^2}= 16r^3\geq 0\quad \mbox{ for every }\, \lambda\in (0, 1],
\end{eqnarray*}
and thus $\frac{\partial\Phi (\lambda, r)}{\partial\lambda}$ is an increasing function of $\lambda$. This gives
\begin{eqnarray*}
\frac{\partial\Phi (\lambda, r)}{\partial\lambda}&\leq &\frac{\partial\Phi }{\partial\lambda}(1, r)= 16r^3-(1-5r)(1-r^2)-8r^2(1+r)\\
&=&-(1-r)^2(1-3r) %\leq 0\quad \mbox{ for }\, r\leq\frac{1}{3}
\end{eqnarray*}
so that $\Phi (\lambda, r)$  is a decreasing function of $\lambda$ on $(0, 1]$ for $r\leq 1/3$ which implies that
\begin{eqnarray*}
\Phi (\lambda, r)&\geq &\Phi (1, r)= %(1-5r)(1-r^2)-8r^2\\ &=&
5r^3-9r^2-5r+1,
\end{eqnarray*}
which is greater than or equal to $0$ for all $r\leq r_*$, where $r_*\approx 0.15867508$ is the unique root of equation
$5r^3-9r^2-5r+1=0,$
which lies in $(0, \frac{1}{5})$.

Since $\Phi (0, r)=2 (1-5r)(1-r^2)$, we have $\Phi (0, r)\geq 0$ for $r\leq 1/5$ and $\Phi (0, r)< 0$ for $r> 1/5$.
Furthermore,
$$\Phi '(1, r)=15r^2-18r-5=15r(r-1)-3r-5<0
$$
which implies that $\Phi (1, r)\geq 0$ for $r\leq r_*$ and $\Phi (1, r)< 0$ for $r> r_*$. According to the fact that $\Phi (\lambda, r)$  is a monotonic decreasing function of $\lambda$ on $(0, 1]$ for $r\leq\frac{1}{3}$, we see that for any $r\in (r_*, \frac{1}{5})$,
$\Phi (0, r)\geq 0,\, \Phi (1, r)<0$, there is a uniquely defined $\lambda(r)\in (0, 1)$ such that $\Phi (\lambda(r), r)=0$.

To prove the last assertion, we have to show that $\frac{d\lambda(r)}{dr}<0$. Indeed, since
$$
\frac{d\lambda(r)}{dr}=-\frac{\frac{\partial\Phi (\lambda , r)}{\partial r}}{\frac{\partial\Phi (\lambda , r)}{\partial\lambda}},
$$
it is sufficient to prove that
\begin{eqnarray*}
\frac{\partial\Phi (\lambda, r)}{\partial r}=24r^2\lambda^2-(39r^2+14r-5)\lambda +30 r^2-4r-10<0
\end{eqnarray*}
for $\lambda\in (0, 1]$ and $r\in (r_*, \frac{1}{5})$.

To that end, we use that for the intervals in question the inequalities
$$\left \{
\begin{array}{l}
30 r^2-4r-10<-8, ~\mbox{and}\\
 24r^2\lambda^2-(39r^2+14r)\lambda+5\lambda-8<24r^2\lambda^2-24r^2\lambda-(15r^2+14r)\lambda <0
\end{array}
\right .
$$
%\begin{eqnarray*}
%30 r^2-4r-10<-8
%\end{eqnarray*}
%and
%\begin{eqnarray*}
%24r^2\lambda^2-(39r^2+14r)\lambda+5\lambda-8<24r^2\lambda^2-24r^2\lambda-(15r^2+14r)\lambda <0
%\end{eqnarray*}
are valid. This completes the proof of Theorem \ref{LSW1-th3}
\epf
% \hfill $\Box$ \vskip 3mm

\bthm\label{LSW1-th4}
%Suppose that $f(z)=h(z)+\overline{g(z)}=\sum_{n=0}^{\infty} a_{n}z^{n}+\overline{\sum_{n=1}^{\infty }b_{n}z^{n}}$ is a sense-preserving harmonic mapping in $\mathbb{D}$ and $h(z)-h(0)\prec_q \varphi(z)-\varphi (0)$ in $\mathbb{D}$, where $\varphi(z)$ is analytic and univalent in $\mathbb{D}$, $\lambda=\dist (\varphi(0), \partial\varphi (\mathbb{D}))<1$ and $\|f_0\|_r=\sum_{n=1}^{\infty}(|a_{n}|^2+|b_{n}|^2)r^{2n}$.
Assume the hypotheses of Theorem {\rm \ref{LSW1-th3}} with a relaxed condition on $\varphi$, namely,  that $\varphi(z)$ is analytic and univalent in $\mathbb{D}$.  Then
\begin{eqnarray*}
T_f(r)=\sum_{n=1}^{\infty}(|a_{n}|+|b_{n}|)r^{n}+\Big(\frac{1}{2-\lambda}+\frac{r}{1-r}\Big)\|f_0\|_r\leq \lambda
\end{eqnarray*}
for $|z|=r\leq r_u^*$, where $r_u^*\approx 0.0808958838$ is the unique root in $(0, 1)$ of the equation
\be\label{LSW1-eq5}
(1-10r+r^2)(1-r)^2(1+r)^3- 32 r^2(1+r^2)=0.
\ee
%The radius $r^*$ is sharp.
\ethm
\bpf %
%Let $\varphi(z)= \sum_{n=0}^{\infty} c_{n}z^{n}$, since $\varphi(z)$ is analytic and univalent in $\mathbb{D}$, by Lemma 1, we have
%$$
%\frac{1}{4}|\varphi'(0)|\leq \lambda\leq |\varphi'(0)|\, \mbox{ and }\, |c_n|\leq n|\varphi'(0)|,\quad n\geq 1,
%$$
%where $\lambda=\dist (\varphi(0), \partial\varphi (\mathbb{D}))$.
Following the notation and the method of proof of Theorem \ref{LSW1-th3}, we easily have $|c_n|\leq 4n\lambda$ for $n\geq 1$
and $T_f(r)\leq 2T_{\varphi}(r)$,
where
\begin{eqnarray*}
2T_{\varphi}(r)&\leq &8\lambda\sum_{n=1}^{\infty}n r^{n}+\Big(\frac{1}{2-\lambda}+\frac{r}{1-r}\Big)32\lambda^2\sum_{n=1}^{\infty}n^2r^{2n}\\
&= &\frac{8\lambda r}{(1-r)^2}+\frac{1+(1-\lambda)r}{(2-\lambda)(1-r)}\cdot \frac{32\lambda^2r^2(1+r^2)}{(1-r^2)^3}\\
&=&\lambda-\lambda\Big[\frac{(1-r)^2-8r}{(1-r)^2}-\frac{32\lambda r^2(1+r^2)(1+(1-\lambda)r)}{(2-\lambda)(1-r)(1-r^2)^3}\Big]\\
&=&\lambda-\lambda\cdot\frac{\Psi (\lambda, r)}{(2-\lambda)(1-r)(1-r^2)^3}.
\end{eqnarray*}
Here
\begin{eqnarray*}
\Psi (\lambda, r)&=&(1-10r+r^2)(2-\lambda)(1-r)^2(1+r)^3- 32\lambda r^2(1+r^2)(1+(1-\lambda)r).
%\\&=& 8r^3\lambda^2-\Big[(1-5r)(1-r^2)+8r^2(1+r)\Big]\lambda +2(1-5r)(1-r^2).
\end{eqnarray*}
To complete the proof, it suffices to  show that $2T_{\varphi}(r)\leq \lambda $ for $r\leq r_u^*$.
%It is clear that $\Psi (\lambda, r)<0$ for $r>1/5$ and $0<\lambda\leq 1$, so that
%$$\lambda-\lambda\cdot\frac{\Psi (\lambda, r)}{(2-\lambda)(1-r)(1-r^2)^3}>\lambda$$for $r>1/5$ and $0<\lambda\leq 1$.

We claim that $\Psi (\lambda, r)\geq 0$ for every $r\leq r_u^*\approx 0.0808958838$ and for $\lambda\in (0, 1]$. In fact, we have
\begin{eqnarray*}
\frac{\partial^2\Psi (\lambda, r)}{\partial\lambda^2}= 64r^3(1+r^2)\geq 0\quad \mbox{ for every }\, \lambda\in (0, 1],
\end{eqnarray*}
and thus $\frac{\partial\Psi (\lambda, r)}{\partial\lambda}$ is an increasing function of $\lambda$. This gives
\begin{eqnarray*}
\frac{\partial\Psi (\lambda, r)}{\partial\lambda}&\leq &\frac{\partial\Psi }{\partial\lambda}(1, r)=-(1-r)\Big[32r^2(1+r^2)+(1-10r+r^2)(1-r)(1+r)^3\Big]\\
%&=&-(1-r)\Big[32r^2(1+r^2)+(1-10r+r^2)(1-r)(1+r)^3\Big]\\
&=&-(1-r)(1-8r+13r^2+51r^4+8r^5-r^6)\leq 0
\end{eqnarray*}
for $r\leq\frac{1}{5}$. Indeed, let $F(r)=1-8r+13r^2+51r^4+8r^5-r^6$. Then we only need to prove $F(r)\geq 0$ for $r\leq\frac{1}{5}$. Since
\begin{eqnarray*}
F'(r)=-8+26r+204r^3+40r^4-6r^5\leq -8+\frac{26}{5}+\frac{204}{125}+\frac{40}{625}<0
\end{eqnarray*}
for $r\leq\frac{1}{5}$, this implies that $F(r)$ is a decreasing function on $[0, 1/5]$ and thus, we conclude that
\begin{eqnarray*}
F(r)&\geq &F(1/5)=1-\frac{8}{5}+\frac{13}{25}+\frac{51}{625}+\frac{8}{5^5}-\frac{1}{5^6}\\
&=&-\frac{2}{25}+\frac{51}{625}+\frac{8}{5^5}-\frac{1}{5^6}=\frac{1}{625}+\frac{39}{5^6}>0
\end{eqnarray*}
for $r\leq\frac{1}{5}$. Hence $\frac{\partial\Psi (\lambda, r)}{\partial\lambda}\leq\frac{\partial\Psi }{\partial\lambda}(1, r)\leq 0$ for $r\leq\frac{1}{5}$. It follows that $\Psi (\lambda, r)$  is a decreasing function of $\lambda$ on $(0, 1]$ for $r\leq\frac{1}{5}$, so that
\begin{eqnarray*}
\Psi (\lambda, r)&\geq &\Psi (1, r)= (1-10r+r^2)(1-r)^2(1+r)^3- 32 r^2(1+r^2),%\\&=&5r^3-9r^2-5r+1,
\end{eqnarray*}
which is greater than or equal to $0$ for all $r\leq r_u^*$, where $r_u^*\approx 0.0808958838$ is the unique root of equation
\eqref{LSW1-eq5}.
%$$(1-10r+r^2)(1-r)^2(1+r)^3- 32 r^2(1+r^2)=0.$$
%The sharpness of $r^*$ can be easily proven by $f(z)=k(z)+\overline{k(z)}$ with $k(z)=z/(1-z)^2$.
This completes the proof of Theorem \ref{LSW1-th4}
\epf

%\hfill $\Box$

\section{The Bohr radius of the derivatives of analytic functions}\label{LPW1-sec5}
%\mbox{}\indent For complex valued analytic function $F(z) =\sum_{n=0}^{\infty} A_n z^n$ defined in $\mathbb{D}$, let the majorant series of $F$ be $M_F (r) := \sum_{n=0}^{\infty} |A_n| r^n,\, |z| = r <1$.
In \cite{BD2019}, Bhowmik and Das investigated the Bohr radius of the derivatives of analytic functions. In particular, they established the following results.

\bprop
{\rm (\cite{BD2019})}
Let $f(z) =\sum_{n=0}^{\infty} a_n z^n$ and $g(z) =\sum_{n=0}^{\infty} b_n z^n$ be two analytic functions in $\mathbb{D}$. Then $M_{f+g}(r)\leq M_{f}(r)+M_{g}(r)$ and $M_{f\, g}(r)\leq M_{f}(r)\, M_{g}(r)$ for any $|z|=r\in [0, 1)$,
where $M_{f}(r)$ denotes the majorant series of $f$.
\eprop

\begin{Thm}\label{LPW1-LemF}
{\rm (\cite{BD2019})}
Let $w(z)$ be an analytic self map of $\mathbb{D}$ with $w(0)=0$. Then $M_{w'}(r)\leq 1$ for $|z| = r \leq r_0 = 1-\sqrt{2/3}$. This radius $r_0$ is the best possible.
\end{Thm}

Using the similar method as in the proof of Theorem \Ref{LPW1-LemF}, we can easily prove the following.

\blem\label{LPW1-lem1}
Let $w(z)$ be an analytic self map of $\mathbb{D}$. Then $M_{zw'}(r)+M_{w}(r)\leq 1$ for $|z| = r \leq r_0 = 1-\sqrt{2/3}$.
\elem

In this section, we determine the Bohr radius for the derivatives of analytic functions associated with quasi-subordination.
More precisely, we have

\bthm\label{LSW1-th5}
Let $f(z) =\sum_{n=0}^{\infty} a_n z^n$ and $g(z) =\sum_{n=0}^{\infty} b_n z^n$ be two analytic functions in $\mathbb{D}$.
If $f(z)-f(0)\prec_q g(z)-g(0)$ in $\mathbb{D}$, then $M_{f'}(r)\leq M_{g'}(r)$ for $|z|=r\leq r_0=1-\sqrt{2/3}$. The radius $r_0$ cannot be improved.
\ethm
\bpf
Suppose that $f(z)-f(0)\prec_q g(z)-g(0)$ in $\mathbb{D}$. Then there exist two functions $\Phi\in {\mathcal B}$ and
$\omega\in {\mathcal B}_0$  such that
$f(z)-f(0) = \Phi(z) (g(\omega(z))-g(0)).
$
Thus we have
$$%%\begin{eqnarray}
f'(z) = \Phi'(z)(g(\omega(z))-g(0))+ \Phi(z) g'(\omega(z)) \omega'(z),
%\label{liu42} \end{eqnarray}
$$
which implies
$$%\begin{eqnarray}
M_{f'}(r)\leq M_{z\Phi'(z)}(r)M_{\frac{\omega(z)}{z}}(r)M_{\frac{g(\omega(z))-g(0)}{\omega(z)}}(r)+M_{\Phi}(r)M_{g'\circ\omega}(r)M_{\omega'}(r).
%\label{liu43} \end{eqnarray}
$$
As $g'\circ\omega\prec g'$ and $\frac{g(\omega(z))-g(0)}{\omega(z)}\prec\frac{g(z)-g(0)}{z}$, by Lemma \Ref{LPW1-LemE}(1), we have $$M_{g'\circ\omega}(r)\leq M_{g'}(r) ~\mbox{ and }~
M_{\frac{g(\omega(z))-g(0)}{\omega(z)}}(r)\leq M_{\frac{g(z)-g(0)}{z}}(r)\leq M_{g'}(r) ~\mbox{ for $r\leq 1/3$}.
$$
From Theorem \Ref{LPW1-LemF}, $M_{\omega'}(r)\leq 1$ for $r\leq r_0=1-\sqrt{2/3}< 1/3$. Further, we observe that
$M_{\frac{\omega(z)}{z}}(r)\leq M_{\omega'}(r)\leq 1$ for $r\leq r_0=1-\sqrt{2/3}< 1/3$. Consequently,
\begin{eqnarray*}
M_{f'}(r)\leq \left (M_{z\Phi'(z)}(r)+M_{\Phi}(r)\right )M_{g'}(r).
%\label{liu44}
\end{eqnarray*}
Moreover Lemma \ref{LPW1-lem1} yields that
\begin{eqnarray*}
M_{z\Phi'(z)}(r)+M_{\Phi}(r)\leq 1 ~ \mbox{ for $r\leq r_0=1-\sqrt{2/3}< 1/3$.} %r M_{\Phi'}(r)+M_{\Phi}(r)\leq 1
%\label{liu45}
\end{eqnarray*}
The desired inequality follows from the last two inequalities.

Following the method of proof of \cite[Theorem 2]{BD2019}, we can easily obtain that the radius $r_0$ cannot be improved.
So, we omit the details. The proof is complete.
\epf
%\hfill $\Box$

\br
It is worth pointing out that  \cite[Theorem 2]{BD2019} is a special case of Theorem \ref{LSW1-th5}.
\er

\subsection*{Acknowledgments}
We thank the referee for his/her careful reading of our paper and invaluable comments. The research of the first author was supported by Guangdong Natural Science Foundation of China (No. 2018A030313508).
The  work of the second author was supported by Mathematical Research Impact Centric Support (MATRICS) of
the Department of Science and Technology (DST), India  (MTR/2017/000367).
The third author was supported by the Natural Science Foundation of China (No. 11771090). The first author would also thank the Laboratory of Mathematics of Nonlinear Sciences, Fudan University (LMNS) for its support during his visit to Fudan University.

%\newpage

\end{document}